\documentclass[11pt,reqno]{amsart}
\usepackage{graphicx, amssymb, color,pdfsync}
\usepackage[all,cmtip]{xy}
\numberwithin{equation}{section}
\usepackage{mathrsfs}
\usepackage{tikz}
\usetikzlibrary{matrix,arrows}

\begin{document}
\newcommand{\s}{\vspace{0.2cm}}

\newtheorem{theo}{Theorem}
\newtheorem{prop}[theo]{Proposition}
\newtheorem{coro}[theo]{Corollary}
\newtheorem{lemm}[theo]{Lemma}
\newtheorem{example}[theo]{Example}
\theoremstyle{remark}
\newtheorem{rema}[theo]{\bf Remark}
\newtheorem{defi}[theo]{\bf Definition}
\newcommand{\rac}{{\mathbb{Q}}}
\newcommand{\comp}{{\mathbb{C}}}
\newcommand{\hip}{{\mathbb{H}}}

\title[The arithmeticity of a Kodaira Fibration]{The arithmeticity of a Kodaira Fibration is determined by its universal cover}
\date{}

\author{Gabino Gonz\'alez-Diez and Sebasti\'an Reyes-Carocca}
\address{Departamento de Matem\'aticas, Universidad Aut\'onoma de Madrid.}
\email{gabino.gonzalez@uam.es, sereyesc@gmail.com}


\thanks{The authors were partially supported by Spanish MEyC Grant MTM 2012-31973. The second author was also partially supported by Becas Chile}

\begin{abstract}
Let $S \to C$ be a Kodaira fibration. Here we show that whether or not the algebraic surface $S$ is defined over a number field depends only on the biholomorphic class of its universal cover. 
\end{abstract}
\maketitle

\section{Introduction and statement of results}
Let $X \subset \mathbb{P}^n$ be a complex  projective variety and $k$ a subfield of the field of the complex numbers $\comp.$ We shall say that $X$ is defined over $k$ or that $k$ is a {\it field of definition} for $X$ if there exists a collection of homogenous polynomials $f_0, \ldots, f_m$ with coefficients in $k$ so that the variety  they define is isomorphic to $X.$ We will say that $X$ is {\it arithmetic} if it is defined over $\overline{\rac}$ or equivalently over a number field.

While it is classically known that there are only three simply connected Riemann surfaces, there is a huge amount of  possibilities for the holomorphic universal cover of a complex surface $S$. It would be interesting to understand the extent to which the arithmeticity of a projective surface can be read off from its holomorphic universal cover. In this short note we study this question for a very important class of complex surfaces known in the literature as Kodaira fibrations.

A {\it Kodaira fibration} consists of a non-singular compact complex surface $S,$ a compact Riemann surface $C$ and a surjective holomorphic map $S \to C$ everywhere of maximal rank such that the fibers are connected and not mutually isomorphic Riemann surfaces. The genera $g$ of the fibre  and $b$ of $C$ are called the genus of the fibration and of the base respectively. It is known that such a surface $S$ must be an algebraic surface of general type and that necessarily $g\ge 3$ and $b \ge 2.$ We notice that an important theorem by Arakelov \cite{Ara} implies that, up to isomorphism, there are only finitely many Kodaira fibrations over a given algebraic curve $C.$

In 1967, Kodaira \cite{Kodaira} used fibrations of this kind to show that the signature of a differentiable fiber bundle need not be multiplicative. Soon after Kas \cite{Kas} studied the deformation space of the surfaces constructed by Kodaira, and two years later Atiyah \cite{Atiyah} and Hirzebruch \cite{Hirzebruch} studied further  properties concerning the signature of Kodaira fibrations in a volume dedicated to Kodaira himself.

Explicit constructions of Kodaira surfaces have been made by Gonz\'alez-Diez and Harvey \cite{HG2}, Bryan and Donagi \cite{BD}, Zaal \cite{Zaal} and Catanese and Rollenske \cite{Catanese}.

We now state the main results of the paper

\begin{theo} \label{theo2} Let $k$ be an algebraically closed  subfield of the complex numbers and $S_1 \to C_1$ and $S_2 \to C_2$ two Kodaira fibrations so that their respective holomorphic universal covers are biholomorphically equivalent. Then $S_1$ is defined over $k$ if and only if $S_2$ is defined over $k.$ In particular, $S_1$ is arithmetic if and only if $S_2$ is arithmetic.
\end{theo}

To prove this theorem we will have to show first the following result which is interesting in its own right

\begin{theo} \label{theo1}  Let $k$ be an algebraically closed  subfield of the complex numbers and $S \to C$ a Kodaira fibration. Then $S$ is defined over $k$ if and only if $C$ is defined over $k.$ In particular, $S$ is an arithmetic surface if and only if $C$ is an arithmetic curve.
\end{theo}

Theorem \ref{theo2} states that the arithmeticity of a Kodaira fibration can be recognized in its holomorphic universal cover. We anticipate that the holomorphic universal cover of $S$ is a contractible bounded domain $\mathscr{B} \subset\comp^2 $ (see Section \ref{s2}). Clearly,
Theorem \ref{theo2}  implies that
  the biholomorphism class of  $\mathscr{B}$ varies together with the variation of $S$ in moduli space. We
  note that in general Kodaira surfaces are not rigid (\cite{Kas}, \cite{Catanese}).

\section{Uniformization of Kodaira Surfaces} \label{s2}

It is well-known that the universal cover of a Riemann surface is biholomorphically equivalent to the projective line $\mathbb{P}^1,$ the complex plane $\comp$ or the upper half-plane $\hip.$ Understanding universal covers of complex manifolds of higher dimension seems to be a very complicated task. However, thanks to the work of Bers \cite{Bers} and Griffiths \cite{Griffiths} on uniformization of algebraic varieties, it is possible to describe the universal cover of a Kodaira fibration $f: S \to C$ in a very explicit way.

Let $\pi : \hip \to C$ be the universal covering map of $C$ and $\Gamma$ the covering group so that $C \cong \hip/ \Gamma.$ By considering the pull-back $h: \pi^*S \to \hip$ of $f$ by $\pi,$ we obtain a new fibration in which, for each $t \in \hip,$ the fiber $h^{-1}(t)$ agrees with the Riemann surface $f^{-1}(\pi(t)).$ Teichm\"{u}ller theory enables us to choose uniformizations $h^{-1}(t)=D_t/K_t$ possessing the following properties:
\begin{itemize}
\item[(a)] $K_t$ is a Kleinian group acting on a bounded domain $D_t$ of $\comp$ which is biholomorphically equivalent to a disk.
\item[(b)] The union of all these disks $\mathscr{B}:=\cup_{t \in \hip} D_t$ is a contractible bounded domain of $\comp^2 $
which is biholomorphic to the universal cover of $S,$ that is, $S \cong \mathscr{B}/\mathbb{G},$ where $\mathbb{G} < \mbox{Aut}(\mathscr{B})$ is the  covering group.
\item[(c)] The group $\mathbb{G}$ is endowed with a surjective homomorphism of groups $\Theta : \mathbb{G} \to \Gamma$ which induces an exact sequence of groups \[ \xymatrix {
  1 \ar[r] &
  \mathbb{K} \ar[r] &
  \mathbb{G} \ar[r]^{\Theta} &
  \Gamma \ar[r] &
  1
} \]
where, for each $t \in \hip$,
 the subgroup $\mathbb{K}$ preserves   $D_t$ and acts on it as the Kleinian group $K_t.$
\end{itemize}

We note that $\mathscr{B}$ carries itself a fibration structure $\mathscr{B} \to \hip$ whose fiber over $t \in \hip$ is $D_t$ (i.e. $\mathscr{B}$ is a {\it Bergman domain} in Bers' terminology).

\

In \cite{Sha} and \cite{Shabat} Shabat studied the automorphism groups of universal covers of families of
 Riemann surfaces and proved a deep result which in
  the case of Kodaira fibrations amounts to the following theorem.

\vspace{0.2cm}
{\bf Theorem (Shabat)}
Let $f: S \to C$ be a Kodaira fibration and $\mathscr{B}$ the holomorphic universal cover of $S.$ Then:
\begin{enumerate}
\item[(a)] the covering group $\mathbb{G}$ of $S$ has finite index in $\mbox{Aut}(\mathscr{B}).$
\item[(b)] $\mbox{Aut}(\mathscr{B})$ is a discrete group.
\end{enumerate}

\section{Proof of Theorems \ref{theo2} and \ref{theo1}} \label{s3}

We denote by $\mbox{Gal}(\comp)$ the group of field automorphisms of $\comp$. The natural action of $\mbox{Gal}(\comp)$ on  the ring of polynomials $\comp[x_0, \ldots, x_n]$ induces a well-defined action $(\sigma, X)\to X^{\sigma}$ on the set of isomorphism classes of algebraic varieties.
Throughout this section $k$ will denote an algebraically closed  subfield of $\comp$ and
  $\mbox{Gal}(\comp/k)$ the subgroup of $\mbox{Gal}(\comp)$ consisting of all automorphisms $\sigma$ fixing the elements of $k.$

\subsection{Proof of Theorem \ref{theo1}} Let $f: S \to C$ be a Kodaira fibration. Let us assume that
the curve $C$ is defined over $k.$ Then  $C^{\sigma}=C$ for all $\sigma \in \mbox{Gal}(\comp/k)$, and so, by Arakelov's finiteness Theorem, there are only finitely many
pairwise non-isomorphic
 Kodaira fibrations $f^{\sigma}: S^{\sigma} \to C^{\sigma}$ with $\sigma \in \mbox{Gal}(\comp/k)$. This implies that $S$ is defined over $k$ \cite[Crit. 2.1]{criterio1}.

\s

In order to prove the converse, we begin by recalling that a complex manifold $X$ is named {\it hyperbolic} if every holomorphic map $\comp \to X$ is a constant map. We claim that  Kodaira fibrations are hyperbolic. In fact, let $f: S \to C$ be a Kodaira fibration and $\varphi : \comp \to S$ a non-constant holomorphic map. As $C$ has genus greater than one, the map $f \circ \varphi: \comp \to C$ must be constant and therefore the image of $\varphi$ has to be contained in a fiber $f^{-1}(x)$ for some $x \in C.$ Since the fibers are also hyperbolic, $\varphi$ must be  constant too.

Let us now assume that $S$ is defined over $k$, so that  $S^{\sigma}=S$ for any $\sigma \in \mbox{Gal}(\mathbb{C}/k)$. Now as $S$ is a
K\"{a}hler hyperbolic  manifold, the canonical divisor $K_S$ is ample \cite{Campana} and this implies that only  finitely many curves  $R$ of genus greater than one can arise as the image of a surjective morphism $S \to R$ \cite{Howard}. In particular the family
 $\{C^{\sigma}:\sigma \in \mbox{Gal}(\mathbb{C}/k)\}$ itself contains only finitely many isomorphism classes of curves. It then follows that $C$ is defined over $k$ \cite[Crit. 2.1]{criterio1}, as required.

\subsection{Proof of Theorem \ref{theo2}} Let $f_2:S_2 \to C_2$ be a Kodaira fibration and $S_1$ an arbitrary non-singular  complex  surface. Let us denote by $\mathscr{B}_i$ the universal cover of $S_i$ and suppose that there exists an isomorphism $\alpha: \mathscr{B}_1 \to \mathscr{B}_2$ between them. Let $\mathbb{G}_i$ be the uniformizing group of $S_i$ so that $\mathscr{B}_i/\mathbb{G}_i \cong S_i.$ By Shabat's Theorem $\mathbb{G}_2$ has finite index in $\mbox{Aut}(\mathscr{B}_2).$ We claim that $\mathbb{G}_1$ has finite index in $\mbox{Aut}(\mathscr{B}_1)$ too. In fact, as $\mathscr{B}_1/\mbox{Aut}(\mathscr{B}_1) \cong \mathscr{B}_2/\mbox{Aut}(\mathscr{B}_2)$ and as $\mbox{Aut}(\mathscr{B}_2)$ is a discrete group, the projection map $S_1=\mathscr{B}_1/\mathbb{G}_1 \to \mathscr{B}_1/\mbox{Aut}(\mathscr{B}_1)$ is a holomorphic map between (normal) compact complex surfaces; from here the claim follows.

By replacing $\mathbb{G}_1$ by $\alpha\mathbb{G}_1\alpha^{-1}$ we can assume that $\mathscr{B}_1=\mathscr{B}_2$,  so we denote $\mathscr{B}_i$ simply by $\mathscr{B}.$ As both $\mathbb{G}_1$ and $\mathbb{G}_2$ have finite index in $\mbox{Aut}(\mathscr{B}),$ so must do their intersection $\mathbb{G}_{12}=\mathbb{G}_{1} \cap \mathbb{G}_{2}.$ The complex surface $S_{12}:=\mathscr{B}/\mathbb{G}_{12}$ is endowed with two finite degree covers $\pi_i':S_{12} \to S_i$ with $i=1,2;$ in particular, $S_{12}$  is also a projective surface. Moreover, if we denote by $\Theta_{12}$ the restriction of the epimorphism $\Theta: \mathbb{G}_2 \to \Gamma_2$ introduced in
   the previous section
  to $\mathbb{G}_{12},$ then we obtain an exact sequence of groups\[ \xymatrix {
  1 \ar[r] &
  \mathbb{K}_{12} \ar[r] &
  \mathbb{G}_{12} \ar[r]^{\Theta_{12}} &
  \Gamma_{12} \ar[r] &
  1
} \]where $\Gamma_{12}=\Theta_{12}(\mathbb{G}_{12})$ and $\mathbb{K}_{12}=\mbox{ker}(\Theta_{12})=\mathbb{K}\cap \mathbb{G}_{12}.$ As in Section \ref{s2}, this sequence defines a Kodaira fibration
$f_{12}: S_{12} \to C_{12}:= \mathbb{H}/\Gamma_{12}$ whose fiber over $[t]_{\Gamma_{12}}$ is  isomorphic to the Riemann surface $D_t/K_t^{12}$ where $K_t^{12}$ is the Kleinian group that realizes the action of $\mathbb{K}_{12}$ on $D_t.$ We have the following commutative diagram
$$
\begin{tikzpicture}[node distance=1.2 cm, auto]
  \node (P) {$\mathscr{B}$};
  \node (A) [below of=P, left of=P] {$S_{12}$};
  \node (C) [below of=B, right of=P] {$S_2$};
    \node (E) [below of=A] {$C_{12}$};
  \node (F) [below of=C] {$C_{2}$};
  \draw[->] (A) to node [swap] {$f_{12}$} (E);
    \draw[->] (P) to node [swap] {} (A);
    \draw[->] (P) to node [swap] {} (C);
        \draw[->] (A) to node {$\pi_2'$} (C);
  \draw[->] (C) to node {$f$} (F);
  \draw[->] (E) to node {$p$} (F);
\end{tikzpicture}
$$where $p$ is the  projection induced by the finite index inclusion $\Gamma_{12} < \Gamma_2.$

Let us now assume that $S_2$ is defined over $k.$ Then Theorem \ref{theo1} ensures that $C_2$ is also defined over $k.$ Furthermore, being an unbranched cover of $C_2,$ the curve $C_{12}$ must also be defined over $k$ \cite[Th. 4.1]{criterio1}. Again, by Theorem \ref{theo1} we  conclude that $S_{12}$ is defined over $k.$ Now, as $S_1$ is a surface of general type arising as the image (by $\pi_1'$) of a surface defined over $k,$ it must be defined over $k$ as well \cite[Prop. 3.2]{criterio1}. This proves Theorem \ref{theo2}.

\end{document}